\documentclass[reqno,11pt]{amsart}
\usepackage{amsmath, amssymb, amsthm,amsfonts} 
\usepackage[english]{babel}
\usepackage{bbm}
\usepackage{graphicx}
\usepackage{url}
\usepackage{epstopdf}
\usepackage{ccaption}
\usepackage{pgfplots}
  \pgfplotsset{compat=newest}
  \usetikzlibrary{plotmarks}
  
\usepackage[a4paper,bindingoffset=0.5cm,left=2cm,right=2cm,top=2.5cm,bottom=2cm,footskip=.8cm]{geometry}

\allowdisplaybreaks

\newcommand{\s}{^\star}

\newcommand{\vb}{\vspace{2.0mm}}

\newcommand{\vt}{\vartheta}

\newtheorem{assumption}{Assumption}
\newtheorem{theorem}{Theorem}
\newtheorem{remark}{Remark}

\makeatletter
\renewcommand{\fnum@figure}[1]{\textbf{\figurename~\thefigure}. }
\renewcommand{\fnum@table}[1]{\textbf{\tablename~\thetable}. }
\makeatother

\begin{document}

\title[An infinite-server system with overdispersed input]{Exact asymptotics in an infinite-server system\\ with overdispersed input}

\author{Mariska Heemskerk \& Michel Mandjes}

\begin{abstract} 
This short communication considers an infinite-server system with overdispersed input. 
The objective is to identify the exact tail asymptotics of the number of customers present at a given point in time under a specific scaling of the model (which involves both the arrival rate and time). 
The proofs rely on a change-of-measure approach. 
The results obtained are illustrated by a series of examples.

\vb

\noindent
{\sc AMS Subject Classification (MSC2010).} 
Primary: 
60K25; 
60F10 
Secondary: 
90B15. 

\vb

\noindent
{\sc Keywords}
Service systems $\circ$ Multi-timescale model $\circ$ Exact asymptotics $\circ$ Overdispersion $\circ$ Staffing.

\vb

\noindent
{\sc Affiliations.} Mariska Heemskerk and Michel Mandjes are with Korteweg-de Vries Institute for Mathematics, University of Amsterdam; Science Park 904, 1098 XH Amsterdam; The Netherlands ({\it email}: {\tt\scriptsize  {j.m.a.heemskerk|m.r.h.mandjes}@uva.nl}). Version: \today.

\vb

\noindent
{\sc Acknowledgments.}   The research of both authors is partly funded by the NWO Gravitation Programme N{\sc etworks} (Grant Number 024.002.003) and an NWO Top Grant (Grant Number 613.001.352). 
 
\end{abstract}

\maketitle

\section{Introduction}
When modeling arrival streams in service systems, the common assumption is that these can be represented by Poisson processes. 
As shown by recent empirical studies, however, this conventional  framework significantly underestimates the level of irregularity present \cite{BAS,LIU,WHI}.  More specifically: arrival streams are often {\it overdispersed}, meaning that the variance of the number of arrivals in a given time window is larger than (rather than equal to) the corresponding mean. 

To remedy this deficiency, various models that lead to overdispersed arrivals have been proposed in the literature. A convenient approach was introduced in \cite{HLM}: every $\Delta>0$ time units, a new value of the arrival rate is {\it sampled} in an i.i.d.\ fashion. More precisely, with $(\Lambda_k)_{k\in{\mathbb N}}$ denoting i.i.d.\ copies of a generic non-negative random variable $\Lambda$, the arrival rate in the interval $I_k:=[(k-1)\Delta,k\Delta)$, for $k\in{\mathbb N}$, is $\Lambda_k.$ 

As infinite-server systems are often used as proxies for their many-server counterparts, they play a prominent role in various staffing rules. This explains the relevance of analyzing the infinite-server system with overdispersed input. Such a model is  studied in e.g.\ \cite{HLM}, where it is assumed that the jobs' service times are  i.i.d.\ samples from some general non-negative distribution (independent of the arrival process), say with distribution function $F(\cdot)$. In \cite{HLM} it is shown that the number of jobs in the system at time $K\Delta$ (assuming the system started empty at time $0$), denoted by $N(K\Delta)$,  is of {\it mixed Poisson} type. More precisely, $N(K\Delta)$ can be written as the sum of $K$ independent terms, i.e., $N_1+\ldots+N_K$, with  $N_k$ having a Poisson distribution with {\it random} parameter 
\[\bar\Lambda_k:=\Lambda_k \int_0^\Delta (1-F((K-k)\Delta+s))\,{\rm d}s;\]
here $\bar\Lambda_k$ corresponds to the contribution of jobs arriving in the interval $I_k$ and still present at time $K\Delta.$ 

Unfortunately, except for some trivial cases there are no closed-form expressions for the distribution of $N(K\Delta)$. A first approach to overcome this is to work with the Laplace transform of $N(K\Delta)$, which turns out to be expressible in terms of the Laplace transform of $\Lambda$, and then to apply numerical inversion (see e.g.\ \cite{AW,dI}). An alternative is to follow a scaling approach, as advocated by \cite{HLM}: after rescaling the random variables $\Lambda_k$ and the sampling interval $\Delta$ in terms of a parameter $n$, explicit characterizations of the distribution of $N(K\Delta)$ can be derived in the asymptotic regime where $n\to\infty$. More specifically, after an appropriate centering and normalization a diffusion limit has been established, as well as rough tail asymptotics (in terms of an exponential decay rate). It is noted, however, that {\it exact} (asymptotically precise, that is) tail asymptotics have not been found so far; such asymptotics can be translated into approximations that can be used when setting up staffing rules for settings in which the desired  service level is relatively high.

The main contribution of this short communication concerns the exact tail asymptotics for the model introduced above; this means that we identify a sequence $f_n$ such that the ratio of the tail probability of interest and $f_n$ converges to $1$ as $n\to\infty$. To be able to work in a convenient framework, we embed the model in a class of (L\'evy-based) models, in a similar way as the procedure followed in \cite{HM}. From that point on, we can {\it mutatis mutandis} follow the line of argumentation that we developed in \cite{HM} to identify the exact asymptotics. Importantly, \cite{HM} focuses on rare events concerning the (overdispersed) arrival process, whereas in this paper we focus on the intrinsically harder counterpart related to the corresponding infinite-server system; indeed the results in this paper are more general than those in \cite{HM}, {in the sense that these can be recovered by sending the service time to $\infty$ (note that a service time almost surely larger than $K\Delta$ suffices)}. 

Importantly, our results also provide a qualitative understanding of the system. For specific choices of our scaling the resampling is so {\it fast} that the infinite-server system effectively experiences Poisson input, so that the asymptotics align with those of a classical M/G/$\infty$ model. In addition there is a region in which the resampling is fast, but not fast enough to provide `full timescale separation'; in that case the asymptotics have to be adapted by adding certain corrections. The opposite regime is the {\it slow} regime, in which the resampling happens relatively infrequently. Again there is the situation in which the timescales {fully} separate, and one in which there is moderate level of timescale separation such that specific corrections appear in the asymptotics. The observed qualitative behavior is in line with the findings in \cite{HM}. 

This short communication is organized as follows. Section \ref{NPS} presents notation, preliminaries and the specific scaling considered; in particular the change of measure featuring in the proofs is introduced. Then, following the setup of \cite{HM}, the fast and slow regime are covered by Sections \ref{FR} and \ref{SR}, respectively. Examples  are presented in Section \ref{NUM}. 

\section{Notation, preliminaries, and scaling} \label{NPS} 
In this section we introduce notation for the infinite-server system described in the introduction. In addition we describe the scaling that we impose throughout this paper, and present the change of measure that will be used in the proofs.

First observe that, using the notation from the introduction and setting $\Lambda(s):=\Lambda_k$ if $s\in I_k$,
\begin{align*}\sum_{k=1}^K\bar\Lambda_k &= \sum_{k=1}^K \Lambda_k \int_0^\Delta (1-F((K-k)\Delta+s))\,{\rm d}s \\&=
\int_0^{K\Delta} \Lambda(s) (1-F(K\Delta-s))\,{\rm d}s\, =_{\rm d}\int_0^{K\Delta} \Lambda(s) (1-F(s))\,{\rm d}s,\end{align*}
with `$=_{\rm d}$' denoting equality in distribution. This representation, in combination with the fact that L\'evy processes can be seen as continuous-time counterparts of random walks, motivates why in this paper we will consider the process 
\[N(t):= A\left(\int_0^{t}  (1-F(s))\,{\rm d}B(s)\right),\]
with $A(\cdot)$ denoting a unit-rate Poisson process and $B(\cdot)$ an increasing L\'evy process (independent of $A(\cdot)$). 
Throughout the paper we work with the {\it characteristic exponents}
\[\alpha(\vt) :=\log {\mathbb E}{\rm e}^{\vt A(1)}={\rm e}^\vt-1,\:\:\:\beta(\vt) :=\log {\mathbb E}{\rm e}^{\vt B(1)},\]
which can be interpreted as  the logarithmic moment generating functions (l-mgf\,s) of $A(1)$ and $B(1)$.
We impose the assumption that $\beta(\cdot)$ is finite in an open neighborhood of the origin, so that we are in a light-tailed regime. The l-mgf of $N(t)$ can be determined by applying standard rules for L\'evy processes, and turns out to equal, with $\bar F(t):=1-F(t)$ denoting the tail distribution 
of the service times,
\[\log{\mathbb E}\,{\rm e}^{\vt N(t)} = \int_0^t \beta\big( \alpha(\vt) \,\bar F(s)\big){\rm d}s.\]
Setting $F(s)\equiv 0$ for all $s\leqslant t$ we recover the l-mgf of $A(B(t))$. 
{As mentioned in the introduction, numerical inversion techniques could in principle help to evaluate the distribution of $N(t)$, but in} this paper we pursue an alternative approach, viz.\ explicit evaluation of the tail probabilities under a specific scaling limit. 

In the scaling limit we consider, time is scaled by a factor $\varphi_n$ and the number of jobs by a factor $n$, as follows. The arrival process we consider is $A(\psi_n B(\varphi_n)\,t)$, assuming that the non-negative sequences $\varphi_n$ and $\psi_n$ are such that $\varphi_n\psi_n=n$ and $\varphi_n\to\infty$ as $n\to\infty$; in the sequel we normalize time such that $t=1$, {which can be done without loss of generality}. The time scaling entails that service times are scaled by $\varphi_n$, such that their distribution function becomes $F(s/\varphi_n).$ The object that we will study is thus
\[N_n:= A\left(\psi_n\int_0^{\varphi_n}  \bar F(s/\varphi_n)\,{\rm d}B(s)\right).\]
The l-mgf $\gamma_n(\cdot)$ of $N_n$ can be expressed in terms of $\alpha(\cdot)$, $\beta(\cdot)$, and $\bar F(\cdot)$:
\[\gamma_n(\vt):=\log {\mathbb E}{\rm e}^{\vt N_n}= \int_0^{\varphi_n} \beta\left(\psi_n \,\alpha(\vt)\,\bar F(s/\varphi_n)\right){\rm d}s=\varphi_n
\int_0^1 \beta\left(\psi_n \,\alpha(\vt)\,\bar F(s)\right){\rm d}s.\]
It requires a straightforward calculation to verify that indeed the number of jobs scales linearly in $n$, in the sense that, with $b:={\mathbb E}\,B(1)=\beta'(0)$,
\[{\mathbb E}\,N_n =\gamma_n'(0) = nc,\:\:\:\mbox{with}\:\:\:c:=b \int_0^1 \bar F(s)\,{\rm d}s.\]
The object of study in this paper is
\[\xi_n(u):= {\mathbb P}(N_n \geqslant un),\]
where we assume that $u>c$ to make sure the event under consideration is {\it rare} (in fact {\it increasingly} rare as $n\to\infty$). More precisely, the focus is on identifying the exact asymptotics of $\xi_n(u)$, meaning that we want to find a sequence $f_n$ such that $\xi_n(u)/f_n\to 1$ as $n\to\infty.$ 

Our analysis is based on a change-of-measure argument. This explains why a crucial role is played by $\vt_n$, defined as the unique positive solution of the equation $\gamma_n'(\vt) = un$; in other words $\vt_n$ solves
\begin{equation}\label{FOC}\int_0^1 \beta'\big(\psi_n\,\alpha(\vt)\,\bar F(s)\big)\,\alpha'(\vt) \,\bar F(s)\,{\rm d}s= u.\end{equation} This $\vt_n$ uniquely exists due to the rarity we assume ($u>c$, that is) in combination with the convexity of $\gamma_n(\cdot)$. The l-mgf of $N_n$ under {the new measure ${\mathbb Q}_n$ can be} expressed in terms 
of the l-mgf of $N_n$ under the original measure, as follows:
\[\gamma^{\mathbb Q_n}_n(\vt) := \gamma_n(\vt+\vt_n)-\gamma_n(\vt_n).\]
This effectively means that twisting $N_n$ by $\vt_n$ leads to a random variable with mean $un$, in the sense that the  measure ${\mathbb Q}_n$ defined through 
\[{\mathbb Q}_n(N_n=k) = {\mathbb P}(N_n=k) {\frac{\exp(\vt_n k)}{{\exp (\gamma_n(\vt_n))}}}\]
has mean $un$; to verify this claim, observe that (by the very definition of $\vt_n$)
\[{\mathbb E}_{{\mathbb Q}_n} N_n=\sum_{k=0}^\infty  k\, {\mathbb Q}_n(N_n= k) = \gamma_n'(\vt_n)=un.\]
For later reference we also compute the variance of $N_n$ under ${\mathbb Q}_n$:
it takes an elementary computation to verify that ${\mathbb V}{\rm ar}_{{\mathbb Q}_n}N_n=\gamma''_n(\vt_n)$ equals 
\begin{equation}
\label{VAR}n\psi_n\int_0^1 \beta''\big(\psi_n\,\alpha(\vt_n)\,\bar F(s)\big)\,(\alpha'(\vt)\bar F(s))^2\,{\rm d}s+n\int_0^1 \beta'\big(\psi_n\,\alpha(\vt_n)\,\bar F(s)\big)\,\alpha''(\vt_n)\,\bar F(s)\,{\rm d}s.\end{equation}

\section{Fast regime}\label{FR}
In this section we consider the case that $\varphi_n$ is superlinear, such that $\psi_n\to 0$ as $n\to\infty.$ This regime is referred to as the {\it fast regime}, as {the timescale corresponding to $B(\cdot)$ is faster than that of the Poisson process $A(\cdot)$}; in the terminology of the introduction, the resampling frequency is relatively high. In our argumentation, we follow the approach developed in \cite[Section 2]{HM}, which borrows elements from the proof of \cite[Thm.\ 3.7.4]{DZ}. The structure of the argumentation is as follows:
\begin{itemize}
\item[$\circ$]
We first analyze the twist factor $\vt_n$, solving $\gamma_n'(\vt) = un.$ As mentioned,  twisting $N_n$ by $\vt_n$ leads to a random variable with mean $un$. It turns out that $\vt_n$ obeys the same type of same expansion as the one featuring in \cite[Section 2]{HM}, i.e.,
\begin{equation}
\vt_n = \sum_{k=0}^\infty v_k\psi_n^{\,k};\label{psif}\end{equation}
evidently, the coefficients $v_k$ are different from those in \cite[Section 2]{HM}, as there only the arrival process was considered (i.e., without jobs potentially leaving the system). 
\item[$\circ$]
The next step is to express  the probability $\xi_n(u)$ using the $\vt_n$-twisted version of $N_n$. By e.g.\ \cite[Ch.\ XIII]{AS},
\begin{equation}\label{COMM}\xi_n(u) ={\mathbb P}(N_n\geqslant un) = {\mathbb E}_{{\mathbb Q}_n}(L(N_n)1\{N_n\geqslant un\}),\end{equation}
with $L(\cdot)$ denoting an appropriate likelihood ratio (translating probabilities under ${\mathbb Q}_n$ into those under the original measure ${\mathbb P}$).  Then the right-hand side of (\ref{COMM}) is further analyzed; from this point on, the proof is identical to that in \cite[Section 2]{HM}.
\end{itemize}
\subsection{Analysis of the twist factor} In this subsection we present a procedure to iteratively find the coefficients $v_k.$ 
The coefficient $v_0$, {which we will refer to as $\vt\s$}, corresponds to $n\to\infty$; using that $\psi_n\to 0$, we find that $\vt\s$ solves
\begin{equation}\label{ths} \beta'(0)\,\alpha'(\vt\s) \,z^+_1 = b\,{\rm e}^{\vt\s} z^+_1= u,\end{equation}
with $z^+_k:=\int_0^1 (\bar F(s))^k\,{\rm d}s.$
We conclude that $\vt\s=\log(u/c)$ (recalling that $c= b\,z^+_1$). Then  $v_1$ can be found from
\[\int_0^1 \beta'(\psi_n\,\alpha(\vt\s+v_1\psi_n)\,\bar F(s))\,\alpha'(\vt\s+v_1\psi_n) \,\bar F(s)\,{\rm d}s= u.\]
Applying Taylor expansions, and using that $\vt\s$ solves (\ref{ths}), we find after some routine calculations that
\begin{equation}\label{v1}v_1= -\frac{\alpha(\vt\s)\alpha'(\vt\s)}{\alpha''(\vt\s)}\frac{\beta''(0)}{\beta'(0)}
\frac{z^+_2}{z^+_1} = -\left(\frac{u}{c}-1\right)\frac{\beta''(0)}{\beta'(0)}
\frac{z^+_2}{z^+_1} .\end{equation}
Using the same ideas, $v_2$ can be expressed in terms of $v_1$. Continuing along the same lines, a procedure can be set up to recursively determine all coefficients $v_k.$
\subsection{Asymptotically exact result}
Equation (\ref{VAR}) reveals that in this fast regime the variance under the new measure ${\mathbb Q}_n$ of $N_n$ grows essentially linearly in $n$, with proportionality constant \[(\sigma_+^{\mathbb Q})^2 := \beta'(0)\alpha''(\vt\s) \,z^+_1 =  b\,{\rm e}^{\vt\s} z^+_1=u .\]
As $\alpha'(\cdot) \equiv \alpha''(\cdot)$, we conclude that under ${\mathbb Q}_n$ the mean and variance of $N_n$ effectively match as $n\to\infty$; cf.\ (\ref{ths}).  This aligns with the heuristic that in the fast regime the resampling is so fast that in essence the system works as an M/G/$\infty$ system (in which the number of jobs has a Poisson distribution); we get back to this intuition below.

\begin{assumption} \label{ass1} 
The sequence $\psi_n$ satisfies
\[\limsup_{n\to\infty} \frac{\log \psi_n}{\log n}<0.\]
\end{assumption}
This assumption entails that there is an $\varepsilon>0$ such that $\psi_n<n^{-\varepsilon}$, and hence $\varphi_n> n^{1+\varepsilon}$, so that  $\varphi_n$ is superlinear. 

We proceed following the argumentation of \cite[Section 2]{HM}; as the line of reasoning is exactly the same, we restrict ourselves to the main steps. The starting point is the  identity
\begin{equation}\label{id}\xi_n(u)= {\mathbb E}_{{\mathbb Q}_n}\left( {\rm e}^{ \gamma_n(\vt_n)-\vt_nN_n}\,1\{N_n \geqslant u n \}\right),\end{equation}
where ${\rm e}^{ \gamma_n(\vt_n)-\vt_n C_n}$ can be interpreted as the likelihood ratio ${\rm d}{\mathbb P}/{\rm d}{{\mathbb Q}_n}$. 
Define 
\[\bar M_n := \frac{N_n -un}{\sqrt{n}\sigma_+^{\mathbb Q}},\]
which has, by the choice of $\vt_n$, mean $0$ under ${\mathbb Q}_n$. Hence, for all $n$,
\begin{equation}\label{ident_fast}
\xi_n(u) = {{\rm e}^{\gamma_n(\vt_n)-\vt_n u n} \, \Delta_n},\:\:\:\mbox{with}\:\:\Delta_n:= {\mathbb E}_{{\mathbb Q}_n}\left( {\rm e}^{-\vt_n\sigma_+^{\mathbb Q}\sqrt{n}\, \bar M_n}\,1\{\bar M_n \geqslant 0 \}\right).\end{equation}
The next step is to  analyze $\delta_n:=\exp(\gamma_n(\vt_n)-\vt_n u n)$ and $\Delta_n$ as $n$ grows large.
\begin{itemize}
\item[$\circ$] First focus on $\delta_n$. Define $m_+\geqslant 1$ through
\[
m_+:=\sup\left\{k\in{\mathbb N}: \liminf_{n\to\infty}\varphi_n\psi_n^{\,k}>0\right\}.
\]
Then we claim, due to  (\ref{psif}), that for appropriately chosen constants ${\bar v}_k$, defining the empty sum as $0$,
\[
\gamma_n(\vt_n)-\vt_n u n=
\chi^+ n +\sum_{k=2}^{m_+} {\bar v}_k \varphi_n\psi_n^{\, k}+o(1),\]
where, recalling that $\alpha(\vt)={\rm e}^\vt-1$, $c =b\,z^+_1$, and $\vt\s=\log(u/c)$,
\begin{equation}\label{chiplus}\chi^+:=b\alpha(\vt\s)\,z^+_1 -\vt\s u=b\big({\rm e}^{\vt\s}-1\big)\,z^+_1 -\vt\s u
=u-c- u\,\log\left(\frac{u}{c}\right)
.\end{equation}
This claim is backed as follows; in passing, the reasoning  shows how the coefficients $\bar v_k$ can be identified. First observe that, expanding $\alpha(\cdot)$ by a Taylor series, $\gamma_n(\vt_n)-\vt_n  un$ equals
\begin{align*}
\varphi_n\int_0^1\beta\left(\psi_n\sum_{\ell=0}^\infty \frac{\alpha^{(\ell)}(\vt\s)}{\ell!}\left(\sum_{k=1}^\infty v_k \psi_n ^{\:k} \right)^\ell \,\bar F(s)\right){\rm d}s-\left(\vt\s +\sum_{k=1}^\infty v_k  \psi_n ^{\:k}\right)un.
\end{align*}
The claim for $m_+=1$  directly follows by expanding $\beta(\cdot)$ through a Taylor series as well,  and collecting terms that are proportional to $n$. For $m_+=2$, ${ \varphi_n\,\psi_n ^{\:k}}\to 0$ when $k> 2$, whereas $\varphi_n\,\psi_n^{\, 2}=n\psi_n$ {stays away from} $0$. As a consequence, including additional terms in the Taylor expansion shows $\gamma_n(\vt_n)-\vt_n  un$ equals, up to terms that are $o(1)$ as $n\to\infty$, 
\[\chi^+n +  \tfrac{1}{2}\beta''(0)(\alpha(\vt\s))^2\,z^+_2\,\varphi_n\psi_n^{\,2}+v_1\left(b\alpha'(\vt\s)\,z^+_1-u\right) \varphi_n\psi_n^{\,2}.\]
This provides us with an expression for $\bar v_2$, using $\alpha(\vt)={\rm e}^\vt-1$ and $\vt\s=\log(u/(bz_1^+))$:
\begin{equation}\label{V2}\bar v_2 = \tfrac{1}{2}\beta''(0)\,(\alpha(\vt\s))^2\,z^+_2+v_1\left(b\alpha'(\vt\s)\,z^+_1-u\right)
=\tfrac{1}{2}\beta''(0)\left(\frac{u}{c}-1\right)^2z_2^+>0.\end{equation}
Higher values of $m_+$ can be dealt with analogously. We have thus developed a procedure to obtain the sequence $(\bar v_k)_{k\geqslant 2}$ from the sequence $(v_k)_{k\geqslant 1}.$ 
\item[$\circ$] We argue how it can be shown  that,  as $n\to\infty$,  $\sqrt{n}\Delta_n$ converges to the positive constant $((1-{\rm e}^{-\vt\s})\sigma_+^{\mathbb Q}\sqrt{2\pi})^{-1}$. 
First, applying integration by parts,
\begin{align}\nonumber
\sqrt{n}\Delta_n &= \sqrt{n}\int_0^\infty {\rm e}^{-\vt_n\sigma_+^{\mathbb Q}\sqrt{n} x}{\mathbb Q}_n(\bar M_n\in{\rm d}x)\\
&={n}\vt_n \sigma_+^{\mathbb Q} \int_0^\infty {\rm e}^{-\vt_n\sigma_+^{\mathbb Q}\sqrt{n} x} \big({\mathbb Q}_n(\bar M_n\leqslant x)-
{\mathbb Q}_n(\bar M_n\leqslant 0)\big){\rm d}x \nonumber \\
&=\sqrt{n}\vt_n\sigma_+^{\mathbb Q}  \int_0^\infty {\rm e}^{-\vt_n\sigma_+^{\mathbb Q} x} \big({\mathbb Q}_n(\bar M_n\leqslant x/\sqrt{n})-
{\mathbb Q}_n(\bar M_n\leqslant 0)\big){\rm d}x.\nonumber
\end{align}
Then it is a matter of applying uniform (in $x$, that is) bounds on ${\mathbb Q}_n(\bar M_n\leqslant x)-\Phi(x)$, with $\Phi(\cdot)$ denoting the cumulative distribution function of a standard Normal random variable; such an Edgeworth expansion is derived in precisely the same way as in \cite[Appendix A]{HM}. Notice that the lattice version, as in \cite[Remark 2]{HM}, needs to be applied, due to the fact that $A(\cdot)$ attains integer values.
\end{itemize}
Combining the above, the following counterpart of \cite[Thm.\ 1]{HM} is obtained.

\begin{theorem}\label{THM1}
As $n\to\infty$, under Assumption \ref{ass1}, 
\[\xi_n(u)\sim \frac{1}{1-{\rm e}^{-\vt\s}}\frac{1}{ \sigma_+^{\mathbb Q} \sqrt{2\pi n}}
\exp\left(\chi^+n +\sum_{k=2}^{m_+} {\bar v}_k \varphi_n\psi_n^{\, k}\right).\]
\end{theorem}

An immediate consequence of Thm.\ \ref{THM1} is that  $\xi_n(u)$ behaves as ${\mathbb P}(A(bz^+_1\,n)\geqslant un)$
when $\varphi_n\psi_n^2 = n\psi_n\to 0$; the process $B(\cdot)$ is effectively replaced by its mean. In this case the exponent is linear in $n$ and equals $\chi^+n$. We are in this situation if for instance  $\varphi_n=n^{f}$ and $f>2$; then the dynamics of $B(\cdot)$ are so much faster than those of $A(\cdot)$  that there  is `full timescale separation'. In addition, Thm.\ \ref{THM1} implies that
the rough (logarithmic) asymptotics are not affected by the choice of $\psi_n$ (as long as Assumption \ref{ass1} is fulfilled): as $n\to\infty$, under Assumption \ref{ass1}, 
\[\frac{1}{n}\log\xi_n(u) \to \chi^+.\]
Observe that $\chi^+$ is the rate function of a Poisson random variable; this once more aligns with the interpretation of {the system in the limit behaving as an M/G/$\infty$ system (whose time-dependent behavior has a Poisson distribution).}

\begin{remark}\label{R1}{\em Whereas in the above reasoning the analysis of $\delta_n$ is relatively straightforward, the analysis of $\Delta_n$ is less intuitive. We therefore include an insightful informal  calculation, based on a discrete version of integration by parts. Let $un$ be integer for simplicity. Write \[\Delta_n = {\rm e}^{\vt_n un}\sum_{k=un}^\infty {\rm e}^{-\vt_n k}{\mathbb Q}_n\left(\bar M_n=\frac{k-un}{\sqrt{n}\sigma_+^{\mathbb Q}}\right).\]
Recall that $\vt_n\to\vt\s$, and observe that ${\rm e}^{-\vt\s k} = \sum_{\ell=k}^\infty {\rm e}^{-\vt\s \ell}(1-{\rm e}^{-\vt\s })$.  Swapping the two sums, we arrive at
\[{\rm e}^{\vt\s un}\sum_{\ell=un}^\infty {\rm e}^{-\vt\s \ell}(1-{\rm e}^{-\vt\s })\sum_{k=un}^\ell {\mathbb Q}_n\left(\bar M_n=\frac{k-un}{\sqrt{n}\sigma_+^{\mathbb Q}}\right).\]
Using the central limit theorem (around the mean; recall that $\bar M_n$ has, under ${\mathbb Q}_n$, mean 0), we approximate
\[\sum_{k=un}^\ell {\mathbb Q}_n\left(\bar M_n=\frac{k-un}{\sqrt{n}\sigma_+^{\mathbb Q}}\right)\approx
(\ell+1-un)\frac{1}{\sqrt{2\pi}}\frac{1}{{\sqrt{n}\sigma_+^{\mathbb Q}}}.\]
Using
${\rm e}^{\vt\s un}\sum_{\ell=un}^\infty {\rm e}^{-\vt\s \ell} (\ell+1-un) = {(1-{\rm e}^{-\vt\s })^{-2}},$
we find $\Delta_n\approx ((1-{\rm e}^{-\vt\s})\sigma_+^{\mathbb Q}\sqrt{2\pi})^{-1}/\sqrt{n}$, as desired.
\hfill$\Diamond$
}\end{remark}

\begin{remark}{\em Above we focus on the exceedance probability $\xi_n(u)$; in this remark we discuss the counterpart of Thm.\ \ref{THM1} that describes the asymptotic behavior of ${\mathbb P}(N_n=un)$ (as in Remark \ref{R1} assuming $un$ is integer). A formal derivation can be given (cf.\ \cite[Exercise  3.7.10]{DZ} with $a=d=1$); we here follow the reasoning of Remark \ref{R1}. The $\delta_n$ is the same as for the exceedance probability case, the counterpart of $\Delta_n$ behaves as
\[{\rm e}^{\vt\s un}\sum_{\ell=un}^\infty {\rm e}^{-\vt\s \ell}(1-{\rm e}^{-\vt\s })\, {\mathbb Q}_n\left(\bar M_n=\frac{k-un}{\sqrt{n}\sigma_+^{\mathbb Q}}\right)\approx
\frac{1}{\sqrt{n}}\frac{1}{\sigma_+^{\mathbb Q}\sqrt{2\pi}};\]
hence the asymptotics differ by a factor $1-{\rm e}^{-\vt\s}$ from those of $\xi_n(\vt)$. Similar properties have been observed in \cite{HKM}.
\hfill$\Diamond$
}\end{remark}

\section{Slow regime}\label{SR}
Here we focus on the case that $\varphi_n$ is sublinear, such that $\psi_n\to \infty$ as $n\to\infty.$ We follow the line of reasoning used in \cite[Section 3]{HM}; observe that (\ref{FOC}) remains valid. The argumentation is as follows:
\begin{itemize}
\item[$\circ$]
In this regime the twist factor $\vt_n$ has the expansion 
\[\vt_n = \sum_{k=1}^\infty w_k\psi_n^{\,-k}.\]
This aligns with the expansion featuring in \cite[Section 3]{HM}, but with different coefficients $w_k$ (that take into account the effect of the leaving jobs). 
\item[$\circ$]
Again, as a next step  the probability $\xi_n(u)$ is rewritten using the $\vt_n$-twisted version of $N_n$; from that point on, the proof precisely follows the one in \cite[Section 3]{HM}.
\end{itemize}
In the sequel we denote
\[z_k^-(\tau) := \int_0^1 \beta^{(k)}\big (\tau \bar F(s)\big) \,(\bar F(s))^k\,{\rm d}s.\]
{Note that the first-order condition $(\ref{FOC})$ can be rewritten as $\mathrm{e}^{\vt}z_1^-(\psi_n(\mathrm{e}^{\vt}-1))=u$.}
\subsection{Analysis of the twist factor} 
Applying Taylor expansions shows that the coefficient $w_1$, which we refer to as $\tau\s$, can be found solving the equation
\[z_1^-(\tau)=\int_0^1 \beta'\big(\tau\,\bar F(s)\big)\,\bar F(s)\,{\rm d}s=u.\]
(where it is used that $\alpha'(0) = 1$). Along the same lines, to identify $w_2$ we find  by first expanding $\alpha(\vt)$ and $\alpha'(\vt)$ in (\ref{FOC}) through Taylor series:
\[\int_0^1 \beta'\left(\big(\tau\s+(w_2+\tfrac{1}{2}(\tau\s)^2)\psi_n^{-1}\big)\,\bar F(s)\right) \,\left(1+\tau\s\psi_n^{-1}\right)\,\bar F(s)\,{\rm d}s +o(\psi_n^{-1})= u\]
(also using that  $\alpha''(0)=1$).
Then expanding $\beta'(\cdot)$ and 
collecting terms of order $\psi_n^{\,-1}$, we obtain
\begin{equation}
\label{W2}w_2=-\tau\s \, \frac{z_1^-(\tau\s)}{z_2^-(\tau\s)}-\tfrac{1}{2}(\tau\s)^2=-\tau\s \, \frac{u}{z_2^-(\tau\s)}-\tfrac{1}{2}(\tau\s)^2>0.\end{equation}
The same procedure can be used to compute the coefficients $(w_k)_{k\geqslant 3}$. 

\subsection{Asymptotically exact result}
By (\ref{VAR}), in this regime the variance under ${\mathbb Q}_n$ of $N_n$ grows essentially linearly in $n\psi_n$, with proportionality constant \[(\sigma_-^{\mathbb Q})^2 := z_2^-(\tau\s)=\int_0^1\beta''(\tau\s \bar F(s))(\bar F(s))^2\,{\rm d}s.\]
The following assumption is the counterpart of Assumption \ref{ass1} for the slow regime.
\begin{assumption}\label{ass2} 
The sequence $\psi_n$ satisfies
\[
0 < \liminf_{n\to\infty} \frac{\log \psi_n}{\log n} \leqslant \limsup_{n\to\infty} \frac{\log \psi_n}{\log n} <1.
\]
\end{assumption}
Due to the first inequality of this assumption there exists an $\varepsilon\in(0,1)$ such that $\psi_n>n^{\varepsilon}$. This implies that $\varphi_n< n^{1-\varepsilon},$ so that $\varphi_n$ is sublinear. 
In addition, by the second inequality also $\psi_n$ is sublinear.

The starting point of the asymptotic analysis of $\xi_n(u)$ is again the identity (\ref{id}). 
We define
\[
\bar E_n := \frac{N_n - u n}{\psi_n\sqrt{\varphi_n}\sigma_-^{\mathbb Q}}
\]
(which has mean $0$ and a variance converging to $1$ under ${\mathbb Q}_n$).
As before, for all $n$, 
\begin{equation}\label{ident_slow}
\xi_n(u) = {\rm e}^{\gamma_n(\vt_n)-\vt_n u n}\, \Delta_n,\:\:\:\mbox{with}\:\:\Delta_n:= {\mathbb E}_{{\mathbb Q}_n}\left( {\rm e}^{-\vt_n\sigma_-^{\mathbb Q}\psi_n\sqrt{\varphi_n}\,\bar E_n}\,1\{\bar E_n \geqslant 0 \}\right).\end{equation}
We are left with analyzing   $\delta_n:=\exp(\gamma_n(\vt_n)-\vt_n u n)$ and  $\Delta_n$ for large $n$.
The analysis of $\delta_n$ can be done as in the fast regime.
Defining
\[
m_-:=\sup\left\{k\in{\mathbb N}: \liminf_{n\to\infty}\varphi_n\psi_n^{-k}>0\right\},
\]
we get, for constants ${\bar w}_k$ and with the empty sum being defined as $0$, as $n\to\infty$,
\[
\delta_n=
\gamma_n(\vt_n)-\vt_n u n=
\chi^-\varphi_n +\sum_{k=1}^{m_-} {\bar w}_k \varphi_n\psi_n^{-k}+o(1),\:\:\:\:\:
\chi^-:=z_0^-(\tau\s)-\tau\s u.
\]
For instance $\bar w_1$ can be identified by collecting the terms that are of order $\varphi_n/\psi_n$; after some algebra this leads to 
\begin{equation}\label{barW1}\bar w_1 =\tfrac{1}{2} (\tau\s)^2 u.\end{equation}
Analogously, the coefficients $(\bar w_k)_{k\geqslant 2}$ can be found.

Following the analysis presented in \cite[Section 3]{HM}, 
$\sqrt{\varphi_n}\Delta_n$ converges to $(\tau\s \sigma_-^{\mathbb Q} \sqrt{2\pi})^{-1}$ as $n\to\infty.$ We thus arrive at the following result, which is the counterpart of \cite[Thm.\ 2]{HM}.

\begin{theorem}\label{THM_SLOW}
As $n\to\infty$, under Assumption \ref{ass2}, for non-lattice $B(\cdot)$,
\[\xi_n(u)\sim \frac{1}{\tau\s \sigma_-^{\mathbb Q} \sqrt{2\pi\varphi_n}}
\exp\left(\chi^-\varphi_n +\sum_{k=1}^{m_-} {\bar w}_k \varphi_n\psi_n^{\, -k}\right).\]
\end{theorem}
From Thm.\ \ref{THM_SLOW} we conclude that if   $\varphi_n\psi_n^{-1}=n/\psi_n^2  \to 0$ as $n\to\infty$, then $\xi_n(u)$ behaves as the probability that the random Poisson parameter
\[ \psi_n \int_0^{\varphi_n} \bar F(s/\varphi_n)\,{\rm d}B(s)\]
exceeds $nu$. If $\varphi_n$ is of the form $n^f$ we have that $\varphi_n\psi_n^{-1} \to 0$ when $f<\frac{1}{2}$. Then the dynamics of $A(\cdot)$ are so much faster than those of $B(\cdot)$ that there is `full timescale separation': the Poisson process is replaced by its rate. In addition, irrespective of the choice of $\psi_n$ (as long as Assumption~\ref{ass2} is met), as $n\to\infty$,
\[\frac{1}{\varphi_n}\log\xi_n(u)\to \chi^-.\]

\section{Examples}\label{NUM}
In this section we present a series of examples illustrating the approximations that we developed. Throughout the process $B(\cdot)$ corresponds to a Gamma process (which is an increasing L\'evy process);  the parameters are $r>0$ (shape) and $\mu>0$ (rate), so that $\beta(\vt) =r  \log\mu-r\log ({\mu - \vt})$ (on the domain $\vt<\mu$).
For the job durations we consider two (crucially distinct) distributions, viz.\ (i)~the exponential distribution, and (ii)~a power-law distribution, and compare with the case that the durations are deterministic. 

First observe that, using that $\alpha(\vt)={\rm e}^\vt-1$,  condition~(\ref{FOC}) can be rewritten as 
\begin{equation}\label{FOC2}
\int_0^1\frac{r\,{\rm e}^\vt\,\bar F(s)}{\mu-\psi_n\,({\rm e}^\vt-1)\,\bar F(s)}\,{\rm d}s=u.
\end{equation}
It takes some algebra to verify that
\begin{equation}\label{zzzz}z_2^-(\tau) = -\frac{z_1^-(\tau)}{\tau} +Z(\tau),\:\:\:\:Z(\tau):=\frac{\mu }{\tau}\int_0^1 \frac{r\,\bar F(s)}{(\mu-\tau\bar F(s))^2}\,{\rm d}s.\end{equation}

\begin{remark}{\em 
Besides a fast and slow regime, there is also the `balanced' regime in which $\varphi_n=n$ and $\psi_n=1.$ It is directly seen that in this case $\gamma_n(\vt)$ is exactly linear in $n$:
\[\gamma_n(\vt) = n \bar\gamma(\vt),\:\:\:\:\bar\gamma(\vt):=\int_0^1 \beta\big(({\rm e}^\vt-1)\bar F(s)\big) \,{\rm d}s=z_0^-({\rm e}^\vt-1).\]
This linearity implies that we are in the framework of the classical Bahadur-Rao result \cite[Thm. 3.7.4]{DZ}. Below we provide the exact asymptotics.
To this end, first observe that in this regime the twist $\vt_n$ does not depend on $n$; it equals $\vt^\circ$, being the solution of
\[\int_0^1 \beta'\big(({\rm e}^\vt-1)\bar F(s)\big) {\rm e}^\vt\bar F(s)\,{\rm d}s = {\rm e}^{\vt^\circ}\,z_1^-({\rm e}^{\vt^\circ}-1)=u.\]
In this case the change-of-measure based derivation of the exact asymptotics is substantially easier than for the fast and slow regime; the argumentation of the proof of \cite[Thm. 3.7.4]{DZ} can be followed. Define
\begin{align}\nonumber(\sigma_\circ^{\mathbb Q})^2 &= \int_0^1 \beta''\big(({\rm e}^{\vt^\circ}-1)\,\bar F(s)\big) ({\rm e}^{\vt^\circ}\bar F(s))^2\,{\rm d}s +
\int_0^1 \beta'\big(({\rm e}^{\vt^\circ}-1)\,\bar F(s)\big) \,{\rm e}^{\vt^\circ}\bar F(s)\,{\rm d}s\\
&={\rm e}^{2\vt^\circ}\,z_2^-({\rm e}^{\vt^\circ}-1)- {\rm e}^{\vt^\circ}\,z_1^-({\rm e}^{\vt^\circ}-1)=
{\rm e}^{2\vt^\circ}\,z_2^-({\rm e}^{\vt^\circ}-1)+u.\label{sigcirc}
\end{align}
We obtain, in line with the findings of \cite[Section 5]{HKM}, with $\chi^\circ:= \bar\gamma(\vt^\circ)-\vt^\circ u$,
\[\xi_n(u)\sim \frac{1}{1-{\rm e}^{-\vt^\circ}}\frac{1}{ \sigma_\circ^{\mathbb Q} \sqrt{2\pi n}}
\exp\left(\chi^\circ n \right).\]
Observe that in this balanced case the exponent is linear in $n$. \hfill$\Diamond$
}\end{remark}

We will consider  the two job-duration distributions mentioned above and show how to determine all the relevant constants for the fast regime with $m_+\in\{1,2\}$ (requiring us to compute the constants $\vt\s$, $\chi^+$, $\bar v_2$, and $(\sigma_{\mathbb Q}^+)^2$), as well as for the balanced regime  (requiring us to compute 
$\vt^\circ$, $\chi^\circ$, and $(\sigma_{\mathbb Q}^\circ)^2$), and for the slow regime with $m_-\in\{0,1\}$  (requiring us to compute
$\tau\s$, $\chi^-$, $\bar w_1$, and $(\sigma_{\mathbb Q}^-)^2$).

\subsection{Exponentially distributed jobs}
We assume exponentially distributed job durations with mean $\nu^{-1}$. To ensure rarity we let $u$ be larger than 
\[c:=\frac{{\mathbb E}N_n}{n}=\frac{\gamma_n'(0)}{n}=\frac{r}{\mu}\frac{1-{\rm e}^{-\nu}}{\nu}.\]
Applying the change-of-variable $t:={\rm e}^{-\nu s}$, some standard calculations yield that (\ref{FOC2}) is equivalent to
\[\frac{1}{\nu}\log\left(\frac{\mu-\psi_n({\rm e}^\vt-1){\rm e}^{-\nu}}{\mu-\psi_n({\rm e}^\vt-1)}\right)=\frac{u\psi_n}{r}(1-{\rm e}^{-\vt}).\]
It is directly verified that $z^+_k= (1-{\rm e}^{-\nu k})/{(\nu k)}$, for $k\in {\mathbb N}$. In addition,
\[z_1^-(\tau) = \frac{r}{\nu\tau} \log\left(\frac{\mu-\tau{\rm e}^{-\nu}}{\mu-\tau}\right),\:\:\:
Z(\tau)=\frac{\mu r}{\nu\tau^2}\left(\frac{1}{\mu-\tau }-\frac{1}{\mu-\tau {\rm e}^{-\nu}}\right).\]
The required constants, pertaining by the three regimes, can be found as follows.

\begin{itemize}
\item[$\circ$]
{\it Fast regime.} 
Recall that in this regime the solution of (\ref{ths}) is given by $\vt\s = \log(u/c)>0.$ We here present the coefficient $\bar v_2$; the next terms in the corresponding expansion can be computed analogously. 
As a consequence of $(\ref{V2})$,
\[\bar v_2 =  \frac{r}{2\mu^2}\left(\frac{u}{c}-1\right)^2\frac{1-{\rm e}^{-2\nu}}{2\nu}.\]
As observed before, in the fast regime $(\sigma_{\mathbb Q}^+)^2=u$; $\chi^+$ is as given by (\ref{chiplus}). 

\item[$\circ$]
{\it Balanced regime.} Here $\vt^\circ$ solves ${\rm e}^{\vt}\,z_1^-({\rm e}^{\vt}-1)=u$, which in this case is equivalent to
\[\frac{1}{\nu}\log\left(\frac{\mu-({\rm e}^\vt-1){\rm e}^{-\nu}}{\mu-({\rm e}^\vt-1)}\right)=\frac{u}{r}(1-{\rm e}^{-\vt});\]
clearly, there is no explicit expression for $\vt^\circ$. 
We do not have an expression for $\chi^\circ$ (in terms of $\vt^\circ$), as (for a given $\vt$) there is no closed form expression for $\bar \gamma(\vt)=z_0^-({\rm e}^\vt-1)$, but it can be expressed in terms of Spence's function ${\rm Li}_2(\cdot)$: by straightforward computations we find
\[\chi^\circ = \frac{r}{\nu}\left({\rm Li}_2\left(\frac{{\rm e}^\vt-1}{\mu}\right)-{\rm Li}_2\left(\frac{({\rm e}^\vt-1){\rm e}^{-\nu}}{\mu}\right)\right)-\vt^\circ u,\]
where ${\rm Li}_2\left(z\right):=-\int_0^z t^{-1}\,\log(1-t)\, \mathrm{d}t$ can be evaluated relying on 
standard numerical techniques. Some calculations, using ${\rm e}^{\vt^\circ}\,z_1^-({\rm e}^{\vt^\circ}-1)=u$, (\ref{zzzz}) and (\ref{sigcirc}), yield an expression for $(\sigma_{\mathbb Q}^\circ)^2$ in terms of $\vt^\circ$:
\begin{equation}
\label{scir}(\sigma_{\mathbb Q}^\circ)^2 = -\frac{u}{{\rm e}^{\vt^\circ}-1}
+
{\rm e}^{2\vt^\circ}\,
Z({\rm e}^{\vt^\circ}-1).\end{equation}
\item[$\circ$]
{\it Slow regime.} With $\tau\s$ solving
$z_1^-(\tau)=u$, by (\ref{barW1}) we  have $\bar w_1 = \tfrac{1}{2} (\tau\s)^2 u$; observe that $\tau\s$ cannot be given explicitly. There is no closed-form expression for $\chi^-$, but again it can be expressed in terms of Spence's function: 
\[\chi^- = \frac{r}{\nu}\left({\rm Li}_2\left(\frac{\tau\s}{\mu}\right)-{\rm Li}_2\left(\frac{\tau\s{\rm e}^{-\nu}}{\mu}\right)\right)-\tau\s u\]
 In addition,
\begin{equation}\label{smin}(\sigma_{\mathbb Q}^-)^2=z_2^-(\tau\s)=-\frac{u}{\tau\s}+Z(\tau\s)
=-\frac{u}{\tau\s}+\frac{\mu r}{\nu(\tau\s)^2}\left(\frac{1}{\mu-\tau\s}-\frac{1}{\mu-\tau\s{\rm e}^{-\nu}  }\right)
.\end{equation}
\end{itemize}

\subsection{Power-law distributed jobs} \label{PLD} We here assume that $\bar F(s) = (1+\kappa s)^{-2}$ for $s\geqslant 0$, with $\kappa>0.$ This distribution is referred to as {\it heavy-tailed} as it has a finite mean $\kappa^{-1}$ but infinite variance. In this case $c=r/(\mu(\kappa+1))$, where we assume that $u>c.$ Define $\eta_n(\vt)\equiv \eta_n := \psi_n({\rm e}^\vt-1)$. A~straightforward computation shows that (\ref{FOC2}) is equivalent to
\[\frac{1}{2\kappa} \sqrt{\frac{\psi_n({\rm e}^\vt-1)}{\mu}} \log\left(\frac{\sqrt{\mu}+\sqrt{\psi_n({\rm e}^\vt-1)}}{\sqrt{\mu}-\sqrt{\psi_n({\rm e}^\vt-1)}}\frac{\sqrt{\mu}(\kappa+1)-\sqrt{\psi_n({\rm e}^\vt-1)}}{\sqrt{\mu}(\kappa+1)+\sqrt{\psi_n({\rm e}^\vt-1)}}\right) = \frac{u\psi_n}{r}(1-{\rm e}^{-\vt}).\]
Also, for $k\in {\mathbb N}$,
\[z_k^+= \frac{1}{2k-1}\frac{1}{\kappa} \left(1-\frac{1}{(1+\kappa)^{2k-1}}\right),\]
whereas (with the derivation of $Z(\cdot)$ in particular taking a considerable amount of calculus)
\begin{align}\nonumber z_1^-(\tau)&=\frac{r}{2\kappa}\frac{1}{\sqrt{\mu\tau}}\log\left(\frac{\sqrt{\mu}+\sqrt{\tau}}{\sqrt{\mu}-\sqrt{\tau}}\frac{\sqrt{\mu}(\kappa+1)-\sqrt{\tau}}{\sqrt{\mu}(\kappa+1)+\sqrt{\tau}}\right),\\ \label{ZP}
Z(\tau) &= \frac{z_1^-(\tau)}{2\tau}+
\frac{r}{2\kappa \tau}\left(\frac{1}{\mu-\tau}- \frac{(\kappa+1)}{\mu(\kappa+1)^2-\tau}\right).
\end{align}
We proceed by considering the fast, balanced, and slow regime.

\begin{itemize}
\item[$\circ$]
{\it Fast regime.}  As before, $\vt\s = \log(u/c)>0$, $(\sigma_{\mathbb Q}^+)^2=u$, and $\chi^+$ is as given by (\ref{chiplus}).
In addition, by $(\ref{V2})$, after some straightforward calculations and using the expression for $z^+_2$,
\[\bar v_2 = \frac{r}{2\mu^2}\left(\frac{u}{c}-1\right)^2\frac{\kappa^2+3\kappa+3}{3(\kappa+1)^3}.\]

\item[$\circ$]
{\it Balanced regime.} As before $\vt^\circ$ solves ${\rm e}^{\vt}\,z_1^-({\rm e}^{\vt}-1)=u$, or equivalently
\begin{equation}\label{balfoc}\frac{1}{2\kappa} \sqrt{\frac{{\rm e}^\vt-1}{\mu}} \log\left(\frac{\sqrt{\mu}+\sqrt{{\rm e}^\vt-1}}{\sqrt{\mu}-\sqrt{{\rm e}^\vt-1}}\frac{\sqrt{\mu}(\kappa+1)-\sqrt{{\rm e}^\vt-1}}{\sqrt{\mu}(\kappa+1)+\sqrt{{\rm e}^\vt-1}}\right) = \frac{u}{r}(1-{\rm e}^{-\vt}).\end{equation}
Again no explicit expressions for $\vt^\circ$ can be given. Relying on (\ref{inte}), and using (\ref{balfoc}),
\begin{align*}\chi^\circ &=     \frac{r}{\kappa}\left((\kappa+1)\log\left(1-\frac{1}{(\kappa+1)^2}\frac{{\rm e}^{\vt^\circ}-1}{\mu}\right)-\log\left(1-\frac{{\rm e}^{{\vt^\circ}}-1}{\mu}\right)\right) \\
&\hspace{1cm}+2({1-{\rm e}^{-\vt^{\circ}}})\,u -\vt^{\circ} u.\end{align*}
The constant 
$(\sigma_{\mathbb Q}^\circ)^2$ can be found by (\ref{scir}), but with $Z(\cdot)$ given by (\ref{ZP}).

\item[$\circ$]
{\it Slow regime.}  Again, by (\ref{barW1}) we  have $\bar w_1 = \tfrac{1}{2} (\tau\s)^2 u$, with  $\tau\s$ solving $z_1^-(\tau)=u$; there is no closed-form expression for $\tau\s$.
Using $z_1^-(\tau\s)=u$ and  (\ref{inte}),
\[\chi^- =\frac{r}{\kappa}\left((\kappa+1)\log\left(1-\frac{1}{(\kappa+1)^2}\frac{\tau\s}{\mu}\right)-\log\left(1-\frac{\tau\s}{\mu}\right)\right) +\tau\s u.\]
The constant 
$(\sigma_{\mathbb Q}^-)^2$  is as in (\ref{smin}), but with $Z(\cdot)$ given by (\ref{ZP}); this leads to
\[(\sigma_{\mathbb Q}^-)^2=-\frac{u}{2\tau\s}+\frac{r}{2\kappa \tau\s}\left(\frac{1}{\mu-\tau\s}- \frac{(\kappa+1)}{\mu(\kappa+1)^2-\tau\s}\right).\]
\end{itemize}

\subsection{Numerical experiments}
In this subsection we report on the numerical experiments carried out for the service-duration distributions discussed above. We evaluate $\xi_n(u)$ for these cases and compare the case of deterministic service times. In the first series of experiments we give the service durations the same mean (namely $\frac{1}{2}$). Note, however, that service durations with the same mean do not necessarily impose the same load on the system; with $c=bz^+_1={\mathbb E}N_n/n$, one could define the load (at time 1) as $c/u$, which we assume to be smaller than 1 to guarantee rarity. To facilitate a comparison under fixed load,  in the second series of experiments we choose the parameters such that for each of the distributions the parameter $z^+_1$ coincides (i.e., $z^+_1=\frac{1}{2}$). 

As mentioned, in the first series of experiments the service times have mean $\frac{1}{2}$, implying that $\nu=\kappa = 2$. 
In these experiments (as well as the ones corresponding to $z^+_1=\frac{1}{2}$)
we present the approximations of $\xi_n(u)$ for  different regimes and levels of timescale separation. More specifically, we present numerical results for $f=\frac{2}{5}$ (slow regime, `full timescale separation' in the sense that $m_-=0$),  $f=\frac{3}{5}$ (slow regime, `moderate timescale separation' in the sense that $m_-=1$),
 $f=1$ (balanced regime),  $f=\frac{5}{3}$ (fast regime, `moderate timescale separation' in the sense that $m_+=2$), and  $f=\frac{5}{2}$ (fast regime, `full timescale separation' in the sense that $m_+=1$). Table \ref{T1} provides the values of all parameters involved in the approximations.

The approximations of $\xi_n(u)$ are given in Table \ref{T2}. For each value of $f$ we chose a corresponding value for $n$ large enough to arrive at tail probabilities roughly of the order $10^{-5}$. 

{\small
\begin{table}[h]
\begin{tabular}{l|cccc|ccc|cccc|}
\cline{2-12}
&$\vt\s$ & $\chi^+$ & $\bar v_2$ & $(\sigma_{\mathbb Q}^+)^2$ &
$\vt^\circ$ & $\chi^\circ$ & $(\sigma_{\mathbb Q}^\circ)^2$ &
$\tau\s$ & $\chi^-$ & $\bar w_1$ & $(\sigma_{\mathbb Q}^-)^2$\\
\hline
Det & 0.693 
& $-0.193$ 
& 0.250
& 1.000
& 0.288 
& $-0.085$ 
& $3.000$
& 0.500
& $-0.153$ 
& 0.125 
& 2.000
\\
Exp & 0.839 
& $-0.271$ 
& 0.212 
& 1.000
& 0.432 
& $-0.150$ 
& 2.608 
& 0.832 
& $-0.319$ 
& 0.346 
& 2.282 
\\
Power-law & 1.099 
& $-0.432$ 
& 0.321 
& 1.000
& 0.551 
& $-0.239$ 
& 3.305 
& 0.961 
& $-0.582$
& 0.461
& 5.977 
\\
\hline
\end{tabular}
\caption{\label{T1}Values of parameters; $r=\mu=u=1$, $\nu=\kappa=2$.}
\end{table}}

{\small
\begin{table}[h] 
\begin{tabular}{l|ccccc|}
\cline{2-6}
$f$ & $\tfrac{2}{5}$ & $\tfrac{3}{5}$ & $1$ & $\tfrac{5}{3}$ & $\tfrac{5}{2}$ \\ \hline
Det 
&$2.613 \cdot 10^{-3}$
&$4.141 \cdot 10^{-3}$
&$1.863 \cdot 10^{-3}$
&$9.644\cdot 10^{-4}$ 
&$4.435\cdot 10^{-4}$
\\
Exp 
&$2.483 \cdot 10^{-5}$	
&$8.188 \cdot 10^{-5}$	
&$5.602 \cdot 10^{-5}$	
&$7.317 \cdot 10^{-5}$	
&$3.792 \cdot 10^{-5} $
\\
Power-law 
& $2.077 \cdot 10^{-8}$
& $1.110 \cdot 10^{-7}$
& $4.693 \cdot 10^{-7}$
& $6.980 \cdot 10^{-7}$
& $2.574 \cdot 10^{-7}$
\\ \hline
$n$ & $3000$ & $200$ & $50$ & $30$ & $30$ \\ 
\hline\end{tabular}
\caption{\label{T2}Approximations of $\xi_n(u)$; mean service time equals $\frac{1}{2}$.}
\end{table}}

However, the table shows that for different service-time distributions with the same mean, the probabilities $\xi_n(u)$ obtained are not necessarily of the same order of magnitude:
the probabilities are highest in the deterministic case ($\gg 10^{-5}$) and, despite its heavy tails, lowest in the power-law case ($\ll 10^{-5}$). 
To explain this ordering, we consider the `loads' corresponding to the three scenarios: observe that $z_1^+ = 0.5$ for deterministic service times, whereas in the exponential case $z_1^+=0.432$, and in the power-law case  $z_1^+ = 0.333$. We thus conclude that the ordering is natural, in the sense that (within each column) the probability $\xi_n(u)$ grows with the system load.

In the second series of experiments, $z^+_1=\frac{1}{2}$ for all service-duration distributions, implying that~$\nu$ is the positive solution of $1-{\rm e}^{-\nu}=\frac{1}{2}\nu$ (so that $\nu\approx 1.594$) and $\kappa=1$; {the deterministic case remains unchanged.} Because we fixed $z_1^+$, the systems have the same load. In Table \ref{T3} the updated values of the parameters are given, while the resulting approximations can be found in Table \ref{T4}. 

{\small
\begin{table}[h]
\begin{tabular}{l|cccc|ccc|cccc|}
\cline{2-12}
&$\vt\s$ & $\chi^+$ & $\bar v_2$ & $(\sigma_{\mathbb Q}^+)^2$ &
$\vt^\circ$ & $\chi^\circ$ & $(\sigma_{\mathbb Q}^\circ)^2$ &
$\tau\s$ & $\chi^-$ & $\bar w_1$ & $(\sigma_{\mathbb Q}^-)^2$\\
\hline
Det & $0.693$ 
& $-0.193$ 
& $0.250$
& $1.000$
& $0.288$ 
& $-0.085$ 
& $3.000$
& $0.500$
& $-0.153$ 
& $0.125$ 
& $2.000$
\\
Exp & $0.693$
& $-0.193$ 
& $0.150$
& $1.000$
& $0.365$
& $-0.108$ 
& $2.358$
& $0.738$
& $-0.236$ 
& $0.272$
& $1.683$ 
\\
Power-law & 0.693
& $-0.193$ 
& 0.146 
& 1.000
& 0.371 
& $-0.109$ 
& 2.313 
& 0.759 
& $-0.243$ 
& 0.288
& 1.668 
\\
\hline
\end{tabular}
\caption{\label{T3}Values of parameters; $r=\mu=u=1$, $\nu=1.594$ and $\kappa = 2$.}
\end{table}}

{\small \begin{table}[h] 
\begin{tabular}{l|ccccc|}
\cline{2-6}
$f$ & $\frac{2}{5}$ & $\frac{3}{5}$ & $1$ & $\frac{5}{3}$ & $\frac{5}{2}$ \\ \hline
Det 
&$3.505 \cdot 10^{-4}$
&$5.304 \cdot 10^{-4}$
&$1.819 \cdot 10^{-4}$
&$4.862\cdot 10^{-5}$ 
&$1.998\cdot 10^{-5}$
\\
Exp 
&$1.296\cdot 10^{-5}$
&$3.193 \cdot 10^{-5}$
&$3.073 \cdot 10^{-5}$
&$3.411 \cdot 10^{-5}$
&$1.998\cdot 10^{-5}$
\\ 
Power-law 
&$9.685 \cdot 10^{-6}$
&$2.516 \cdot 10^{-5}$
&$2.732 \cdot 10^{-5}$
&$3.356 \cdot 10^{-5}$
&$1.998 \cdot 10^{-5}$
\\ \hline
$n$ & $8000$ & $400$ & $75$ & $45$ & $45$ \\ 
\hline\end{tabular}
\caption{\label{T4}Approximations of $\xi_n(u)$; $z^+_1$ equals $\frac{1}{2}$.}
\end{table}}


Indeed, this table  shows that for the three service-time distributions the probabilities $\xi_n(u)$ are of roughly the same order of magnitude. Observe that in this setting we chose different (larger) values for $n$ than before, to again guarantee probabilities roughly of the order $10^{-5}$ (note that a larger mean results in a higher probability of exceeding level $ u n$).  
In this setting with constant load, one would have perhaps anticipated that $\xi_n(u)$ is largest in the power-law case (due to its heavy tail) and smallest in the deterministic case. 
Realize however that in the time domain considered (i.e., $[0,1]$) the tails of the distributions do not play a significant role yet; we refer to \cite{HKM} for related findings. 

\appendix
\section{Auxiliary computations for power-law distribution}
In this appendix we concentrate on computing $z_0^-(\tau)$ for the case of our power-law distributed service times featuring in Section \ref{PLD}. Denoting $F_\tau:= (r/2\kappa)\,\sqrt{\tau/\mu}$ and $T_\tau(x):=\tau/(\mu(x+1)^2)$, by a change-of-variables argument, 
\[z_0^-(\tau) = - F_\tau\int_{T_\tau(\kappa)}^{T_\tau(0)} \frac{\log(1-t)}{t\sqrt{t}}\,{\rm d}t.\]
Applying integration by parts (using that the primitive of $t^{-3/2}$ is $-2\,t^{-1/2}$), this integral equals
\[ \frac{r}{\kappa}\left(\log\left(1-\frac{\tau}{\mu}\right)-(\kappa+1)\log\left(1-\frac{1}{(\kappa+1)^2}\frac{\tau}{\mu}\right) \right)+2F_\tau\int_{T_\tau(\kappa)}^{T_\tau(0)} \frac{1}{1-t}\frac{1}{\sqrt{t}}\,{\rm d}t.\]
Using the identity
\[ \int_{T_\tau(\kappa)}^{T_\tau(0)} \frac{1}{1-t}\frac{1}{\sqrt{t}}\,{\rm d}t = \log\left(\frac{\sqrt{\mu}+\sqrt{\tau}}{\sqrt{\mu}-\sqrt{\tau}}\frac{\sqrt{\mu}(\kappa+1)-\sqrt{\tau}}{\sqrt{\mu}(\kappa+1)+\sqrt{\tau}}\right)=\frac{2\kappa}{r}{\sqrt{\mu\tau}} \,z_1^-(\tau), \]
we conclude that $z_0^-(\tau) $ equals 
\begin{equation}
\label{inte}
  \frac{r}{\kappa}\left(\log\left(1-\frac{\tau}{\mu}\right)-(\kappa+1)\log\left(1-\frac{1}{(\kappa+1)^2}\frac{\tau}{\mu}\right)\right) +2\tau \,z_1^-(\tau).
\end{equation}

\end{document}